# THE BEHAVIOR OF THE NPMLE OF A DECREASING DENSITY NEAR THE BOUNDARIES OF THE SUPPORT

By Vladimir N. Kulikov and Hendrik P. Lopuhaä

*Eurandom and Delft University of Technology*

We investigate the behavior of the nonparametric maximum likelihood estimator $\hat{f}_n$ for a decreasing density $f$ near the boundaries of the support of $f$. We establish the limiting distribution of $\hat{f}_n(n^{-\alpha})$, where we need to distinguish between different values of $0 < \alpha < 1$. Similar results are obtained for the upper endpoint of the support, in the case it is finite. This yields consistent estimators for the values of $f$ at the boundaries of the support. The limit distribution of these estimators is established and their performance is compared with the penalized maximum likelihood estimator.

**1. Introduction.** In various statistical models, such as density estimation and estimation of regression curves or hazard rates, monotonicity constraints can arise naturally. For these situations certain isotonic estimators have been in use for considerable time. Often these estimators can be seen as maximum likelihood estimators in a semiparametric setting. Although conceptually these estimators have great appeal and are easy to formulate, their distributional properties are usually of a very complicated nature.

In the context of density estimation, the nonparametric maximum likelihood estimator $\hat{f}_n$ for a nonincreasing density $f$ on $[0, \infty)$ was studied by Grenander [2]. It is defined as the left derivative of the least concave majorant (LCM) of the empirical distribution function $F_n$ constructed from a sample from $f$. Prakasa Rao [11] obtained the asymptotic pointwise behavior of $\hat{f}_n$. Groeneboom [3] provided an elegant proof of the same result, which can be formulated as follows. For each $x_0 > 0$,

$$(1.1) \quad |4f(x_0)f'(x_0)|^{-1/3} n^{1/3}\{\hat{f}_n(x_0) - f(x_0)\} \to \underset{t \in \mathbb{R}}{\arg\max}\{W(t) - t^2\}$$

in distribution, where $W$ denotes standard two-sided Brownian motion originating from zero. The first distributional result for a global measure of deviation for $\hat{f}_n$ was found by Groeneboom [3], concerning asymptotic normality of the $L_1$-distance $\|\hat{f}_n - f\|_1$ (see [4] for a rigorous proof).

Apart from estimating a monotone density $f$ on $(0, \infty)$, the estimation of the value of $f$ or its derivatives at zero is required in various statistical applications. There is a direct connection with renewal processes, where the backward recurrence time in equilibrium has density $f(x) = (1 - G(x))/\mu$, where $G$ and $\mu$ are the distribution function and mean of the interarrival times (see [1]). Clearly, $f$ is decreasing and a natural parameter of interest is $\mu = 1/f(0)$. An interesting application is in the context of natural fecundity of human populations, where one is interested in the time $T$ it takes for a couple from initiating attempts to become pregnant until conception occurs. Keiding, Kvist, Hartvig and Tvede [6] investigated a current-duration design where data are collected from a cross-sectional sample of couples that are currently attempting to become pregnant. If $U$ is the time to discontinuation without success and $V$ is the time to discontinuation of follow-up, then $X = T \wedge U$ is the waiting time until termination for whatever reason, and $Y = T \wedge U \wedge V$ is the observed experience waiting time. When the initiations happen according to a homogeneous Poisson process, $Y$ is distributed as the backward recurrence time in a renewal process in equilibrium, and the survival function of $X$ is $f(x)/f(0)$, where $f$ is decreasing. Woodroofe and Sun [13] provide a different application in the context of astronomy. If $Y$ denotes the normalized angular diameter of a galaxy, conditional on that it is being observed, then $1/Y^3$ has a nonincreasing density $f$ and the proportion of galaxies that are observed is $1/f(0)$. Another example is from Hampel [5], who studied the sojourn time of migrating birds. Under certain model assumptions, the expected sojourn time is $-f(0)/f'(0)$, where $f$ is the (convex) decreasing density of the time span between capture and recapture of a bird.

In contrast to (1.1), Woodroofe and Sun [13] showed that $\hat{f}_n$ is not consistent at zero. They proposed a penalized maximum likelihood estimator $\hat{f}_n^{\mathrm{P}}(0)$ and in [12] it was shown that

$$n^{1/3}\{\hat{f}_n^{\mathrm{P}}(0) - f(0)\} \to \sup_{t>0} \frac{W(t) - (c - f(0)f'(0)t^2/2)}{t},$$

where $c$ depends on the penalization. Surprisingly, the inconsistency of $\hat{f}_n$ at zero does not influence the behavior of $\|\hat{f}_n - f\|_1$. Nevertheless, the inconsistency at the boundaries will have an effect if one studies other global measures of deviation, such as the $L_k$-distance, for $k$ larger than 1, or the supremum distance.

In this paper we study the behavior of the Grenander estimator at the boundaries of the support of $f$. We first consider a nonincreasing density $f$ on $[0, \infty)$ and investigate the behavior of

$$(1.2) \qquad n^{\beta}\{\hat{f}_n(cn^{-\alpha}) - f(cn^{-\alpha})\}$$

for $c > 0$, where $0 < \alpha < 1$ and $\beta > 0$ are chosen suitably in order to make (1.2) converge in distribution. Our results will imply that when $f'(0) < 0$, then $\hat{f}_n(cn^{-1/3})$ is a consistent estimator for $f(0)$ at rate $n^{1/3}$ with a limiting distribution that is a functional of $W$. This immediately yields $\hat{f}_n^{\mathrm{S}}(0) = \hat{f}_n(n^{-1/3})$ as a simple estimator for $f(0)$. A more adaptive alternative would be to find the value of $c$ that minimizes the asymptotic mean squared error. This turns out to depend on $f$ and then has to be estimated. The resulting estimator $\hat{f}_n^{\mathrm{A}}(0) = \hat{f}_n(\hat{c}n^{-1/3})$ will be compared with the penalized maximum likelihood estimator from [12]. We will also consider the case where $f'(0) = 0$ and $f''(0) < 0$, which requires different values for $c$ and $\alpha$. For nonincreasing $f$ with compact support, say $[0, 1]$, we also investigate the behavior near 1. Similarly, this leads to a consistent estimator for $f(1)$. Moreover, the results on the behavior of $\hat{f}_n$ at the boundaries of $[0, 1]$ allow an adequate treatment of the $L_k$-distance between $\hat{f}_n$ and $f$. It turns out that for $k > 2.5$, the inconsistency of $\hat{f}_n$ starts to affect the behavior of $\|\hat{f}_n - f\|_k$ (see [10]).

In Section 2 we give a brief outline of our approach for studying differences such as (1.2) and state some preliminary results for the arg max functional. Section 3 is devoted to the behavior of $\hat{f}_n$ near zero. Section 4 deals with the behavior of $\hat{f}_n$ near the boundary at the other end of the support for a density $f$ on $[0, 1]$. In Section 5 we propose two estimators $\hat{f}_n^{\mathrm{S}}(0)$ and $\hat{f}_n^{\mathrm{A}}(0)$ based on the presented theory, and compare these with the penalized maximum likelihood estimator from Sun and Woodroofe [12].

**2. Preliminaries.** Instead of studying the process $\{\hat{f}_n(t) : t \geq 0\}$ itself, we will use the more tractable inverse process $\{U_n(a) : a \geq 0\}$, where $U_n(a)$ is defined as the last time that the process $F_n(t) - at$ attains its maximum,

$$U_n(a) = \underset{t \in [0, \infty)}{\arg\max}\{F_n(t) - at\}.$$

Its relation with $\hat{f}_n$ is as follows: with probability 1

$$(2.1) \qquad \hat{f}_n(x) \leq a \quad \Longleftrightarrow \quad U_n(a) \leq x.$$

Let us first describe the line of reasoning used to prove convergence in distribution of (1.2). We illustrate things for the case $c = 1$, $0 < \alpha < 1/3$, and $f'(0) < 0$. It turns out that in this case the proper choice for $\beta$ is $1/3$. Hence, we will consider events of the type

$$n^{1/3}\{\hat{f}_n(n^{-\alpha}) - f(n^{-\alpha})\} \leq x.$$

According to relation (2.1), this event is equivalent to

$$U_n(f(n^{-\alpha}) + xn^{-1/3}) - n^{-\alpha} \leq 0.$$

The left-hand side is the arg max of the process

$$Z_n(t) = F_n(t + n^{-\alpha}) - f(n^{-\alpha})t - xtn^{-1/3}.$$

With suitable scaling, the process $Z_n$ converges in distribution to some Gaussian process $Z$. The next step is to use an arg max version of the continuous mapping theorem from [7]. The version that suffices for our purposes is stated below for further reference.

THEOREM 2.1. *Let $\{Z(t) : t \in \mathbb{R}\}$ be a continuous random process satisfying:*

  (i) *$Z$ has a unique maximum with probability 1.*
  (ii) *$Z(t) \to -\infty$, as $|t| \to \infty$, with probability 1.*

*Let $\{Z_n(t) : t \in \mathbb{R}\}$ be a sequence of random processes satisfying:*

  (iii) *$\arg\max_{t \in \mathbb{R}} Z_n(t) = O_p(1)$, as $n \to \infty$.*

*If $Z_n$ converges in distribution to the process $Z$ in the topology of uniform convergence on compacta, then $\arg\max_{t \in \mathbb{R}} Z_n(t)$ converges in distribution to $\arg\max_{t \in \mathbb{R}} Z(t)$.*

This theorem yields that $U_n(f(n^{-\alpha}) + xn^{-1/3})$, properly scaled, converges in distribution to the arg max of a Gaussian process. Convergence of (1.2) then follows from another application of (2.1).

The main difficulty in verifying the conditions of Theorem 2.1 is showing that (iii) holds. It requires careful handling of all small order terms in the expansion of the process. In the process of proving condition (iii) we will frequently use the following lemma, which enables us to suitably bound the arg max from above.

LEMMA 2.1. *Let $f$ and $g$ be continuous functions on $K \subset \mathbb{R}$.*

  (i) *Suppose that $g$ is nonincreasing. Then $\arg\max_{x \in K}\{f(x) + g(x)\} \leq \arg\max_{x \in K} f(x)$.*
  (ii) *Let $C > 0$ and suppose that for all $s, t \in K$, such that $t \geq C + s$, we have that $g(t) \leq g(s)$. Then $\arg\max_{x \in K}\{f(x) + g(x)\} \leq C + \arg\max_{x \in K} f(x)$.*

In studying processes like $Z_n$ we will use a Brownian approximation similar to the one used in [4]. Let $E_n$ denote the empirical process $\sqrt{n}(F_n - F)$. For $n \geq 1$, let $B_n$ be versions of the Brownian bridge constructed on the same probability space as the uniform empirical process $E_n \circ F^{-1}$ via the Hungarian embedding, where

$$(2.2) \qquad \sup_{t \geq 0} |E_n(t) - B_n(F(t))| = O_p(n^{-1/2} \log n)$$

(see [8]). Define versions $W_n$ of Brownian motion by

$$W_n(t) = B_n(t) + \xi_n t, \qquad t \in [0, 1],$$

where $\xi_n$ is a standard normal random variable independent of $B_n$. This means that we can represent $B_n$ by the pathwise equality $B_n(t) = W_n(t) - tW_n(1)$.

We will often apply a Brownian scaling argument in connection with $\arg\max$ functionals. Note that $\arg\max_t\{Z(t)\}$ does not change by multiplying $Z$ by a constant, and that the process $W(bt)$ has the same distribution as the process $b^{1/2}W(t)$. This implies that

$$
\begin{aligned}
a \arg\max_{t \in I}\{W(bt) - ct^k\} &= \arg\max_{t \in aI}\{W(ba^{-1}t) - ca^{-k}t^k\} \\
(2.3) \qquad &\stackrel{d}{=} \arg\max_{t \in aI}\{b^{1/2}a^{-1/2}W(t) - ca^{-k}t^k\} \\
&= \arg\max_{t \in aI}\{W(t) - cb^{-1/2}a^{-k+1/2}t^k\}
\end{aligned}
$$

for $I \subset \mathbb{R}$ and constants $a, b > 0$ and $c \in \mathbb{R}$.

**3. Behavior near zero.** We first consider the case that $f$ is a nonincreasing density on $[0, \infty)$ satisfying:

(C1) $0 < f(0) = \lim_{x \downarrow 0} f(x) < \infty$.

(C2) For some $k \geq 1$, $0 < |f^{(k)}(0)| \leq \sup_{s \geq 0} |f^{(k)}(s)| < \infty$, with $f^{(k)}(0) = \lim_{x \downarrow 0} f^{(k)}(x)$, and $f^{(i)}(0) = 0$ for $1 \leq i \leq k - 1$.

Under these conditions we determine the behavior of the Grenander estimator near zero. With the proper normalizing constants the limit distribution of $n^\beta(\hat{f}_n(n^{-\alpha}) - f(n^{-\alpha}))$ is independent of $f$. Define $D[Z(t)](a)$ as the right derivative of the LCM on $\mathbb{R}$ of the process $Z(t)$ at the point $t = a$, and define $D_R$ similarly, where the LCM is restricted to the set $t \geq 0$.

THEOREM 3.1. *Suppose $f$ satisfies conditions* (C1) *and* (C2) *and let $c > 0$. Then:*

(i) *For $1/(2k + 1) < \alpha < 1$ and $A_1 = (c/f(0))^{1/2}$, the sequence*

$$A_1 n^{(1-\alpha)/2}(\hat{f}_n(cn^{-\alpha}) - f(cn^{-\alpha}))$$

*converges in distribution to $D_R[W(t)](1)$ as $n \to \infty$.*

(ii) *For* $\alpha = 1/(2k+1)$, $B_{2k} = (f(0)^{1/2}|f^{(k)}(0)|^{-1}(k+1)!)^{2/(2k+1)}$ *and* $A_{2k} = \sqrt{B_{2k}/f(0)}$, *the sequence*

$$A_{2k}\left\{ n^{(1-\alpha)/2}(\hat{f}_n(cB_{2k}n^{-\alpha}) - f(cB_{2k}n^{-\alpha})) + \frac{f^{(k)}(0)(cB_{2k})^k}{k!} \right\}$$

*converges in distribution to* $D_{\mathrm{R}}[W(t) - t^{k+1}](c)$ *as* $n \to \infty$.

(iii) *For* $0 < \alpha < 1/(2k+1)$ *and* $A_{3k} = (2(k-1)!)^{1/3}|f(0)f^{(k)}(0)c^{k-1}|^{-1/3}$, *the sequence*

$$A_{3k}n^{1/3+\alpha(k-1)/3}(\hat{f}_n(cn^{-\alpha}) - f(cn^{-\alpha}))$$

*converges in distribution to* $D[W(t) - t^2](0)$ *as* $n \to \infty$.

REMARK 3.1. In order to present the limiting distributions in Theorem 3.1 in the same way, they have been expressed in terms of slopes of least concave majorants. However, note that similar to the switching relation (2.1), one finds that

$$D_{\mathrm{R}}[W(t)](1) \stackrel{d}{=} \sqrt{\underset{t \in [0,\infty)}{\arg\max}\{W(t) - t\}},$$

$$D[W(t) - t^2](0) \stackrel{d}{=} 2\underset{t \in \mathbb{R}}{\arg\max}\{W(t) - t^2\}.$$

In studying the behavior of (1.2), we follow the line of reasoning described in Section 2. We start by establishing convergence in distribution of the relevant processes. It turns out that we have to distinguish between three cases concerning the rate at which $n^{-\alpha}$ tends to zero.

LEMMA 3.1. *Suppose* $f$ *satisfies* (C1) *and* (C2) *and let* $W$ *denote standard two-sided Brownian motion on* $\mathbb{R}$. *For* $1/(2k+1) \le \alpha < 1$, $t \ge 0$ *and* $x \in \mathbb{R}$, *define*

$$Z_{n1}(x,t) = n^{(1+\alpha)/2}(F_n(tn^{-\alpha}) - f(0)tn^{-\alpha}) - xt.$$

(i) *For* $1/(2k+1) < \alpha < 1$, *the process* $\{Z_{n1}(x,t) : t \in [0,\infty)\}$ *converges in distribution in the uniform topology on compacta to the process* $\{W(f(0)t) - xt : t \in [0,\infty)\}$.

(ii) *For* $\alpha = 1/(2k+1)$, *the process* $\{Z_{n1}(x,t) : t \in [0,\infty)\}$ *converges in distribution in the uniform topology on compacta to* $\{W(f(0)t) - xt + f^{(k)}(0)t^{k+1}/(k+1)! : t \in [0,\infty)\}$.

(iii) *For* $0 < \alpha < 1/(2k+1)$, $b = (1 - 2\alpha(k-1))/3$, $t \ge -cn^{b-\alpha}$ *and* $x \in \mathbb{R}$, *define*

$$Z_{n2}(x,t) = n^{(b+1)/2}(F_n(cn^{-\alpha} + tn^{-b}) - F_n(cn^{-\alpha}) - f(cn^{-\alpha})tn^{-b}) - xt.$$

*Then the process $\{Z_{n2}(x,t) : t \in [-cn^{b-\alpha}, \infty)\}$ converges in distribution in the uniform topology on compacta to the process $\{W(f(0)t) - xt + c^{k-1}f^{(k)}(0)t^2/ (2(k-1)!) : t \in \mathbb{R}\}$.*

The next step is to use Theorem 2.1. The major difficulty is to verify condition (iii) of this theorem. The following lemma ensures that this condition is satisfied.

LEMMA 3.2. *Let $f$ satisfy* (C1) *and* (C2) *and let $Z_{n1}$, $Z_{n2}$ and $b$ be defined as in Lemma* 3.1.

(i) *For $1/(2k+1) < \alpha < 1$ and $x > 0$, $\arg\max_{t \in [0,\infty)} Z_{n1}(x,t) = \mathcal{O}_p(1)$.*

(ii) *For $\alpha = 1/(2k+1)$ and $x \in \mathbb{R}$, $\arg\max_{t \in [0,\infty)} Z_{n1}(x,t) = \mathcal{O}_p(1)$.*

(iii) *For $0 < \alpha < 1/(2k+1)$ and $x \in \mathbb{R}$, $\arg\max_{t \in [-cn^{b-\alpha}, \infty)} Z_{n2}(x,t) = \mathcal{O}_p(1)$.*

With Lemmas 3.1 and 3.2 at hand, the proof of Theorem 3.1 consists of using the switching relation (2.1) and an application of Theorem 2.1.

PROOF OF THEOREM 3.1. (i) First note that by condition (C2),

$$n^{(1-\alpha)/2}(\hat{f}_n(cn^{-\alpha}) - f(cn^{-\alpha})) = n^{(1-\alpha)/2}(\hat{f}_n(cn^{-\alpha}) - f(0))$$
$$+ \mathcal{O}(n^{(1-(2k+1)\alpha)/2}),$$

where $(1 - (2k+1)\alpha)/2 < 0$. For $x > 0$, according to (2.1),

$$(3.1) \qquad \begin{aligned} P\{n^{(1-\alpha)/2}(\hat{f}_n(cn^{-\alpha}) - f(0)) \le x\} \\ = P\{n^\alpha U_n(f(0) + xn^{-(1-\alpha)/2}) \le c\}. \end{aligned}$$

If $Z_{n1}$ is the process defined in Lemma 3.2(i), then

$$(3.2) \qquad 0 \le n^\alpha U_n(f(0) + xn^{-(1-\alpha)/2}) = \arg\max_{t \in [0,\infty)} Z_{n1}(x,t) = \mathcal{O}_p(1),$$

where, according to Lemma 3.1, the process $\{Z_{n1}(x,t) : t \in [0,\infty)\}$ converges in distribution to the process $\{W(f(0)t) - xt : t \in [0,\infty)\}$. To apply Theorem 2.1, we have to extend the above processes to the whole real line. Therefore define

$$\tilde{Z}_{n1}(t) = \begin{cases} Z_{n1}(x,t), & t \ge 0, \\ t, & t \le 0. \end{cases}$$

Then for $x$ fixed, $\tilde{Z}_{n1}$ converges in distribution to the process $Z_1$, where

$$Z_1(t) = \begin{cases} W(f(0)t) - xt, & t \ge 0, \\ t, & t \le 0. \end{cases}$$

Moreover, since $Z_{n1}(x, 0) = 0$, together with (3.2), it follows that

$$\underset{t \in \mathbb{R}}{\arg \max} \, \tilde{Z}_{n1}(t) = \underset{t \in [0, \infty)}{\arg \max} \, \tilde{Z}_{n1}(t)$$

$$= n^\alpha U_n(f(0) + xtn^{-(1-\alpha)/2}) = \mathcal{O}_p(1).$$

The process $Z_1$ is continuous, and since $\operatorname{Var}(Z_1(s) - Z_1(t)) \neq 0$ for $s, t > 0$ with $s \neq t$, it follows from Lemma 2.6 in [7] that $Z_1$ has a unique maximum with probability 1. By an application of the law of the iterated logarithm for Brownian motion,

$$(3.3) \qquad P\left\{ \limsup_{|u| \to \infty} \frac{W(u)}{\sqrt{2|u| \log \log |u|}} = 1 \right\} = 1,$$

it can be seen that $Z_1(t) \to -\infty$ as $|t| \to \infty$. Theorem 2.1 now yields that $\arg \max_{t \in \mathbb{R}} \tilde{Z}_{n1}(t)$ converges in distribution to

$$\underset{t \in \mathbb{R}}{\arg \max} \, Z_1(t) = \underset{t \geq 0}{\arg \max} \{ W(f(0)t) - xt \}.$$

Using (3.1) together with (2.3), this implies that

$$P\{ n^{(1-\alpha)/2}(\hat{f}_n(cn^{-\alpha}) - f(0)) \leq x \}$$

$$= P\left\{ \underset{t \in \mathbb{R}}{\arg \max} \, \tilde{Z}_{n1}(t) \leq c \right\}$$

$$\to P\left\{ \underset{t \geq 0}{\arg \max} \{ W(f(0)t) - xt \} \leq c \right\}$$

$$= P\left\{ \underset{t \geq 0}{\arg \max} \left\{ W(t) - \frac{xc^{1/2}t}{f(0)^{1/2}} \right\} \leq 1 \right\}.$$

Similar to the switching relation (2.1), the right-hand side equals

$$P\{ (f(0)/c)^{1/2} D_{\mathrm{R}}[W(t)](1) \leq x \},$$

so that it remains to show that $P\{ n^{(1-\alpha)/2}(\hat{f}_n(cn^{-\alpha}) - f(0)) \leq 0 \} \to 0$. But this is evident, as for any $\varepsilon > 0$, using (2.3) once more,

$$P\{ n^{(1-\alpha)/2}(\hat{f}_n(cn^{-\alpha}) - f(0)) \leq 0 \}$$

$$\leq P\{ n^{(1-\alpha)/2}(\hat{f}_n(cn^{-\alpha}) - f(0)) \leq \varepsilon \}$$

$$\to P\left\{ \underset{t \geq 0}{\arg \max} \left\{ W(t) - \frac{\varepsilon t}{\sqrt{f(0)}} \right\} \leq c \right\}$$

$$= P\left\{ \underset{t \geq 0}{\arg \max} \{ W(t) - t \} \leq \frac{c\varepsilon^2}{f(0)} \right\}.$$

When $\varepsilon \downarrow 0$, the right-hand side tends to zero, which can be seen from

$$P\left\{\limsup_{t \downarrow 0} \frac{W(t)}{\sqrt{2t \log\log(1/t)}} = 1\right\} = 1.$$

This proves (i).

(ii) First note that by (C2),

$$n^{k/(2k+1)}(\hat{f}_n(cB_{2k}n^{-1/(2k+1)}) - f(cB_{2k}n^{-1/(2k+1)})) + f^{(k)}(0)\frac{(cB_{2k})^k}{k!}$$

$$= n^{k/(2k+1)}(\hat{f}_n(cB_{2k}n^{-1/(2k+1)}) - f(0)) + o(1),$$

and that according to (2.1), $P\{n^{k/(2k+1)}(\hat{f}_n(cB_{2k}n^{-1/(2k+1)}) - f(0)) \le x\}$ is equal to

$$P\{B_{2k}^{-1}n^{1/(2k+1)}U_n(f(0) + xn^{-k/(2k+1)}) \le c\}.$$

With $Z_{n1}$ being the process defined in Lemma 3.1 with $\alpha = 1/(2k+1)$, we get

$$B_{2k}^{-1}n^{1/(2k+1)}U_n(f(0) + xn^{-k/(2k+1)}) = \operatorname*{arg\,max}_{t \in [0,\infty)}\{Z_{n1}(x, B_{2k}t)\} = \mathcal{O}_p(1).$$

Again we first extend the above process to the whole real line:

$$\tilde{Z}_{n1}(t) = \begin{cases} Z_{n1}(x, B_{2k}t), & t \ge 0, \\ t, & t \le 0. \end{cases}$$

Then, according to Lemma 3.1, $\tilde{Z}_{n1}$ converges in distribution to the process

$$Z_2(t) = \begin{cases} W(f(0)B_{2k}t) - B_{2k}xt + f^{(k)}(0)B_{2k}^{k+1}t^{k+1}/(k+1)!, & t \ge 0, \\ t, & t \le 0. \end{cases}$$

Similar to the proof of (i), it follows from Theorem 2.1 that $\operatorname{arg\,max}_t \tilde{Z}_{n1}(t)$ converges in distribution to $\operatorname{arg\,max}_t Z_2(t)$. This implies that

$$P\{A_{2k}n^{k/(2k+1)}(\hat{f}_n(cB_{2k}n^{-1/(2k+1)}) - f(0)) \le x\}$$

$$\to P\left\{\operatorname*{arg\,max}_{t \ge 0}\left\{W(f(0)B_{2k}t) - \frac{B_{2k}xt}{A_{2k}} + \frac{f^{(k)}(0)B_{2k}^{k+1}t^{k+1}}{(k+1)!}\right\} \le c\right\}$$

$$= P\left\{\operatorname*{arg\,max}_{t \ge 0}\{W(t) - xt - t^{k+1}\} \le c\right\}$$

$$= P\{D_{\mathrm{R}}[W(t) - t^{k+1}](c) \le x\},$$

by means of Brownian scaling similar to (2.3), and a switching relation similar to (2.1).

(iii) According to (2.1), we have

$$
(3.4)
\begin{aligned}
&P\{n^{(1-b)/2}(\hat{f}_n(cn^{-\alpha}) - f(cn^{-\alpha})) \le x\} \\
&\qquad = P\{n^b(U_n(f(cn^{-\alpha}) + xn^{-(1-b)/2}) - cn^{-\alpha}) \le 0\},
\end{aligned}
$$

and with $Z_{n2}$ as defined in Lemma 3.2(iii), we get

$$
n^b(U_n(f(cn^{-\alpha}) + xn^{-(1-b)/2}) - cn^{-\alpha}) = \operatorname*{arg\,max}_{t \in [-cn^{b-\alpha}, \infty)} Z_{n2}(x, t) = \mathcal{O}_p(1).
$$

As in the proof of (i) and (ii), we extend the above process to the whole real line:

$$
\tilde{Z}_{n2}(t) = \begin{cases} Z_{n2}(x, t), & t \ge -cn^{b-\alpha}, \\ Z_{n2}(x, -cn^{b-\alpha}) + (t + cn^{b-\alpha}), & t < -cn^{b-\alpha}. \end{cases}
$$

Then by Lemma 3.1 $Z_{n2}$ converges in distribution to the process $Z_3$, where

$$
Z_3(t) = W(f(0)t) - xt + \frac{f^{(k)}(0)c^{k-1}}{2(k-1)!}t^2, \qquad t \in \mathbb{R}.
$$

Similar to the proofs of (i) and (ii), it follows from Theorem 2.1 that $\arg\max_t Z_{n2}(t)$ converges in distribution to $\arg\max_t Z_3(t)$. Together with (3.4), this implies that

$$
\begin{aligned}
&P\{n^{(1-b)/2}A_{3k}(\hat{f}_n(cn^{-\alpha}) - f(cn^{-\alpha})) \le x\} \\
&\qquad \to P\left\{\operatorname*{arg\,max}_{t \in \mathbb{R}}\left\{W(f(0)t) - A_{3k}^{-1}xt + \frac{f^{(k)}(0)c^{k-1}}{2(k-1)!}t^2\right\} \le 0\right\} \\
&\qquad = P\left\{\operatorname*{arg\,max}_{t \in \mathbb{R}}\{W(t) - xt - t^2\} \le 0\right\} \\
&\qquad = P\{D[W(t) - t^2](0) \le x\},
\end{aligned}
$$

again using Brownian scaling similar to (2.3), and a switching relation similar to (2.1).  □

## 4. Behavior near the end of the support.

Suppose that $f$ has compact support and, without loss of generality, assume this to be the interval $[0, 1]$. In this section we investigate the behavior of $\hat{f}_n$ near 1. Although there seems to be no simple symmetry argument to derive the behavior near 1 from the results in Section 3, the arguments to obtain the behavior of

$$
n^\beta\{f(1 - n^{-\alpha}) - \hat{f}_n(1 - n^{-\alpha})\}
$$

are similar to the ones used in studying (1.2). If $f(1) > 0$, then $\hat{f}_n(1)$ will always underestimate $f(1)$, since by definition $\hat{f}_n(1) = 0$. Nevertheless, the behavior near the end of the support is similar to the behavior near zero.

For this reason, we only provide the statement of a theorem for the end of the support, which is analogous to Theorem 3.1. For details on the proof we refer to [9]. Motivations for studying the behavior near the end of the support are not so strong as for the behavior near zero. However, the behavior near 1 is required for establishing the asymptotic normality of the $L_k$-distance between $\hat{f}_n$ and $f$. Similar to (C1) and (C2) we will assume that:

(C3)  $0 < f(1) = \lim_{x \uparrow 1} f(x) < \infty$.

(C4)  For some $k \geq 1$, $0 < |f^{(k)}(1)| \leq \sup_{0 \leq s \leq 1} |f^{(k)}(s)| < \infty$, with $f^{(k)}(1) = \lim_{x \uparrow 1} f^{(k)}(x)$ and $f^{(i)}(1) = 0$ for $1 \leq i \leq k - 1$.

We then have the following theorem.

THEOREM 4.1.   *Suppose $f$ satisfies conditions* (C3) *and* (C4) *and* $c > 0$. *Then:*

(i)  *For $1/(2k+1) < \alpha < 1$ and $\tilde{A}_1 = (c/f(1))^{1/2}$, the sequence*

$$\tilde{A}_1 n^{(1-\alpha)/2} (f(1 - cn^{-\alpha}) - \hat{f}_n(1 - cn^{-\alpha}))$$

*converges in distribution to $D_{\mathrm{R}}[W(t)](1)$ as $n \to \infty$.*

(ii)  *For $\alpha = 1/(2k+1)$, $\tilde{B}_{2k} = (f(1)^{1/2}|f^{(k)}(1)|^{-1}((k+1)!))^{2/(2k+1)}$ and $\tilde{A}_{2k} = \sqrt{\tilde{B}_{2k}/f(1)}$, the sequence*

$$\tilde{A}_{2k} \left\{ n^{(1-\alpha)/2} (f(1 - c\tilde{B}_{2k} n^{-\alpha}) - \hat{f}_n(1 - c\tilde{B}_{2k} n^{-\alpha})) - \frac{|f^{(k)}(1)|(c\tilde{B}_{2k})^k}{k!} \right\}$$

*converges in distribution to $D_{\mathrm{R}}[W(t) - t^{k+1}](c)$ as $n \to \infty$.*

(iii)  *For $0 < \alpha < 1/(2k+1)$ and $\tilde{A}_{3k} = ((k-1)!)^{1/3}|4f(1)f^{(k)}(1)c^{k-1}|^{-1/3}$, the sequence*

$$\tilde{A}_{3k} n^{1/3 + \alpha(k-1)/3} (f(1 - cn^{-\alpha}) - \hat{f}_n(1 - cn^{-\alpha}))$$

*converges in distribution to $D[W(t) - t^2](0)$ as $n \to \infty$.*

PROOF.   The proof is similar to that of Theorem 3.1. We briefly sketch the proof for case (i); details can be found in [9].

Similar to the proof of Theorem 3.1(i), it suffices to consider

$$n^{(1-\alpha)/2} (f(1) - \hat{f}_n(1 - cn^{-\alpha})).$$

For $x > 0$, according to (2.1),

$$(4.1) \qquad \begin{aligned} &P\{n^{(1-\alpha)/2}(f(1) - \hat{f}_n(1 - cn^{-\alpha})) \leq x\} \\ &\quad = P\{n^\alpha (1 - U_n(f(1) - xn^{-(1-\alpha)/2})) \leq c\}. \end{aligned}$$

We have that $n^\alpha(1 - U_n(f(1) - xn^{-(1-\alpha)/2})) = \arg\max_{t \in [0,n^\alpha]} Y_{n1}(x,t)$, where the process

$$Y_{n1}(x,t) = n^{(1+\alpha)/2}(F_n(1 - tn^{-\alpha}) - F_n(1) + f(1)tn^{-\alpha}) - xt$$

converges in distribution to the process $\{W(f(1)t) - xt : t \in [0,\infty)\}$. From here on, the proof proceeds in completely the same manner as that of Theorem 3.1(i). We conclude that for $x > 0$,

$$P\{n^{(1-\alpha)/2}(f(1) - \hat{f}_n(1 - cn^{-\alpha})) \le x\}$$
$$= P\Big\{\arg\max_{0 \le t \le n^\alpha} Y_{n1}(t) \le c\Big\}$$
$$\to P\Big\{\arg\max_{t \ge 0}\{W(f(1)t) - xt\} \le c\Big\}$$
$$= P\Big\{\arg\max_{t \ge 0}\Big\{W(t) - \frac{xc^{1/2}t}{f(1)^{1/2}}\Big\} \le 1\Big\}.$$

By (2.1), the right-hand side equals $P\{(f(1)/c)^{1/2}D_R[W(t)](1) \le x\}$, and similar to the proof of Theorem 3.1(i) it follows that $P\{n^{(1-\alpha)/2}(f(1) - \hat{f}_n(1 - cn^{-\alpha})) \le 0\} \to 0$. This proves (i). □

## 5. A comparison with the penalized NPMLE.

Consider a decreasing density $f$ on $[0,\infty)$. We first consider the case where $f'(0) < 0$. As pointed out in [13], the NPMLE $\hat{f}_n$ for $f$ is not consistent at zero. They proposed a penalized NPMLE $\hat{f}_n^P(\alpha_n, 0)$, and in Sun and Woodroofe [12] they show that

$$n^{1/3}\{\hat{f}_n^P(\alpha_n, 0) - f(0)\} \to \sup_{t > 0} \frac{W(t) - (c - (1/2)f(0)f'(0)t^2)}{t},$$

where $c$ is related to the smoothing parameter $\alpha_n = cn^{-2/3}$. Sun and Woodroofe [12] also provide (to some extent) an adaptive choice for $c$ that leads to an estimate $\hat{\alpha}_n$ of the smoothing parameter, and report some results of a simulation experiment for $\hat{f}_n^P(\hat{\alpha}_n, 0)$.

We propose two consistent estimators of $f(0)$, both converging at rate $n^{1/3}$. A simple estimator is $\hat{f}_n^S(0) = \hat{f}_n(n^{-1/3})$. This estimator is straightforward and does not have any additional smoothing parameters. According to Theorem 3.1(ii), $\hat{f}_n^S(0)$ is a consistent estimator for $f(0)$, converging at rate $n^{1/3}$. It has a limiting distribution that is a functional of $W$,

$$A_{21}n^{1/3}\{\hat{f}_n^S(0) - f(0)\} \to D_R[W(t) - t^2](1/B_{21}),$$

where $A_{21}$ and $B_{21}$ are defined in Theorem 3.1(ii). In order to reduce the mean squared error, we also propose an adaptive estimator

$$\hat{f}_n^A(0) = \hat{f}_n(c_1^* \hat{B}_{21} n^{-1/3})$$

for $f(0)$. Here $c_k^*$ is the value that minimizes $E(D_R[W(t) - t^{k+1}](c))^2$, and $\hat{B}_{21}$ is an estimate for the constant $B_{21}$ in Theorem 3.1(ii). Computer simulations show that $c_k^* \approx 0.345$ for both $k = 1$ and $k = 2$. We take

$$\hat{B}_{21} = 4^{1/3} \hat{f}_n^S(0)^{1/3} |\tilde{f}_n'(0)|^{-2/3},$$

where

$$\tilde{f}_n'(0) = \min(n^{1/6}(\hat{f}_n(n^{-1/6}) - \hat{f}_n(n^{-1/3})), -n^{-1/3})$$

is an estimate for $f'(0)$. As we have seen above, $\hat{f}_n^S(0)$ is consistent for $f(0)$, and according to Theorem 3.1, $\tilde{f}_n'(0)$ is consistent for $f'(0)$. When $f$ is twice continuously differentiable, it converges at rate $n^{1/6}$. Therefore $\hat{B}_{21}$ is consistent for $B_{21}$ and $\hat{f}_n^A(0)$ is a consistent estimator of $f(0)$, converging with rate $n^{1/3}$. It has the limit behavior

$$A_{21} n^{1/3} \{\hat{f}_n^A(0) - f(0)\} \to D_R[W(t) - t^2](c_1^*),$$

where $A_{21}$ is defined in Theorem 3.1(ii).

We simulated 10,000 samples of sizes $n = 50$, 100, 200 and 10,000 from a standard exponential distribution with mean 1. For each sample, the values of $n^{1/3}\{\hat{f}_n^S(0) - f(0)\}$, $n^{1/3}\{\hat{f}_n^A(0) - f(0)\}$ and $n^{1/3}\{\hat{f}_n^P(\hat{\alpha}_n, 0) - f(0)\}$ were computed. The value of $\hat{\alpha}_n$ was computed as proposed in [12], $\hat{\alpha}_n = 0.649 \cdot \hat{\beta}_n^{-1/3} n^{-2/3}$, where

$$\hat{\beta}_n = \max\left\{\hat{f}_n^P(\alpha_0, 0) \frac{\hat{f}_n^P(\alpha_0, 0) - \hat{f}_n^P(\alpha_0, x_m)}{2x_m}, n^{-q}\right\}$$

is an estimate of $\beta = -f(0)f'(0)/2$. Here $x_m$ denotes the second point of jump of the penalized NPMLE $\hat{f}_n^P(\alpha_0, \cdot)$ computed with smoothing parameter $\alpha_0$. The parameter $\alpha_0 = c_0 n^{-2/3}$, and $q$ should be taken between 0 and 0.5. However, Sun and Woodroofe [12] do not specify how to choose $q$ and $c_0$ in general. We took $q = 1/3$, and for $\alpha_0$ the values as listed in their Table 2: $\alpha_0 = 0.0516$, 0.0325 and 0.0205 for sample sizes $n = 50$, 100 and 200. For sample size $n = 10,000$ we took the theoretical optimal value $\alpha_0 = 0.649\beta^{-1/3} n^{-2/3}$, with $\beta = 0.5$. It is worth noticing that Sun and Woodroofe [12] do not optimize the MSE, but $n^{1/3} E|\hat{f}_n^P(\hat{\alpha}_n, 0) - f(0)|$. Nevertheless, computer simulations show that the $\alpha_n$ minimizing the MSE is approximately the same and that $n^{2/3} E|\hat{f}_n^P(\alpha, 0) - f(0)|^2$ is a very flat function in a neighborhood of $\alpha_n$. A similar property holds for the value $c_k^*$ minimizing the AMSE of our estimator.

In Table 1 we list simulated values for the mean, variance and mean squared error of the three estimators. The penalized NPMLE is less biased, but has a larger variance. Estimator $\hat{f}_n^A(0)$ performs better in the sense of mean squared error, approaching the best theoretically expected performance. It is also remarkable how well it mimics its limiting distribution for small samples. Estimator $\hat{f}_n^{S}(0)$ performs a little worse than $\hat{f}_n^{A}(0)$, having the largest bias, but the smallest variance.

If $k = 2$ in condition (C2), it is possible to estimate $f(0)$ at a rate faster than $n^{1/3}$. If it is known in advance that $k = 2$, we can produce two consistent estimators of $f(0)$ converging at rate $n^{2/5}$. Similar to the previous case, a simple estimator is $\hat{f}_n^{S,2}(0) = \hat{f}_n(n^{-1/5})$. It is a consistent estimator of $f(0)$, converging at rate $n^{2/5}$, and has the limit behavior

$$A_{22} n^{2/5} \{\hat{f}_n^{S,2}(0) - f(0)\} \to D_R[W(t) - t^3](1/B_{22}),$$

where $A_{22}$ and $B_{22}$ are defined in Theorem 3.1(ii). Again, we propose an adaptive estimator $\hat{f}_n^{A,2}(0) = \hat{f}_n(c_2^* \hat{B}_{22} n^{-1/5})$ for $f(0)$, where $\hat{B}_{22}$ is an estimate for the constant $B_{22} = 36^{1/5} f(0)^{1/5} |f''(0)|^{-2/5}$ in Theorem 3.1(ii), and $c_2^* \approx 0.345$ is the value that minimizes $E(D_R[W(t) - t^3](c))^2$. We take $\hat{B}_{22} = 36^{1/5} \hat{f}_n^{S,2}(0)^{1/5} |\tilde{f}_n''(0)|^{-2/5}$, where we estimate $f''(0)$ by $\tilde{f}_n''(0) = \min(2n^{1/4} \times (\hat{f}_n(n^{-1/8}) - \hat{f}_n(n^{-1/5})), -n^{-1/8})$. As we have seen above, $\hat{f}_n^{S,2}(0)$ is consistent for $f(0)$, and according to Theorem 3.1, $\tilde{f}_n''(0)$ is consistent for $f''(0)$ with rate $n^{1/8}$ if $f$ is three times continuously differentiable. Therefore $\hat{B}_{22}$ is a consistent estimator for $B_{22}$ and $\hat{f}_n^{A,2}(0)$ is a consistent estimator of $f(0)$, converging with rate $n^{2/5}$:

$$A_{22} n^{2/5} \{\hat{f}_n^{A,2}(0) - f(0)\} \to D_R[W(t) - t^3](c_2^*),$$

where $A_{22}$ is defined in Theorem 3.1(ii).

We simulated 10,000 samples of sizes $n = 50$, 100, 200 and 10,000 from a half-normal distribution. For each sample, the values of $n^{2/5} \{\hat{f}_n^{S,2}(0) - f(0)\}$ and $n^{2/5} \{\hat{f}_n^{A,2}(0) - f(0)\}$ were computed. Sun and Woodroofe [12] do not consider the possibility of constructing a special estimator for the case $k = 2$, though we believe that this is also possible with a penalization technique. In Table 2 we list simulated values for the mean, variance and mean squared error of both estimators. The simple estimator is more biased but its variance is smaller than the variance of the adaptive one.

TABLE 1
*Simulated mean, variance and mean squared error for the three estimators at the standard exponential distribution*

|  |  | $n$ | | | |
| --- | --- | --- | --- | --- | --- |
|  |  | 50 | 100 | 200 | 10,000 |
| $n^{1/3}\{\hat{f}_n^{S}(0) - f(0)\}$ | Mean | −0.847 | −0.853 | −0.868 | −0.917 |
|  | Var | 0.439 | 0.484 | 0.536 | 0.700 |
|  | MSE | 1.157 | 1.211 | 1.289 | 1.541 |
| $n^{1/3}\{\hat{f}_n^{A}(0) - f(0)\}$ | Mean | −0.738 | −0.777 | −0.793 | −0.643 |
|  | Var | 0.934 | 0.742 | 0.807 | 1.045 |
|  | MSE | 1.478 | 1.345 | 1.436 | 1.458 |
| $n^{1/3}\{\hat{f}_n^{P}(\hat{\alpha}_n, 0) - f(0)\}$ | Mean | −0.072 | −0.079 | −0.075 | −0.195 |
|  | Var | 1.296 | 1.530 | 1.732 | 1.913 |
|  | MSE | 1.301 | 1.537 | 1.738 | 1.951 |

If it is not known in advance that $k = 2$, then application of estimators $\hat{f}_n^{\mathrm{S},2}(0)$ and $\hat{f}_n^{\mathrm{A},2}(0)$ is undesirable. If in fact $k = 1$, they are still consistent, but their convergence rate will be $n^{1/5}$. On the other hand, when $k = 2$, then $\hat{f}_n^{\mathrm{S}}(0)$, $\hat{f}_n^{\mathrm{A}}(0)$ and $f_n^{\mathrm{P}}(\hat{\alpha}_n, 0)$ are still applicable. In that case, according to Theorem 3.1(i), $\hat{f}_n^{\mathrm{S}}(0)$ is a consistent estimator of $f(0)$ converging at rate $n^{1/3}$, such that

$$n^{1/3}\{\hat{f}_n^{\mathrm{S}}(0) - f(0)\} \to \sqrt{f(0)} D_{\mathrm{R}}[W(t)](1).$$

Also $\hat{f}_n^{\mathrm{A}}(0)$ is still consistent for $f(0)$ in case $k = 2$, but now at rate $n^{7/18}$. This can be seen as follows. Since $f'(0) = 0$, it follows that

$$n^{1/6}\tilde{f}_n'(0) \to -\sqrt{f(0)} D_{\mathrm{R}}[W(t)](1) + \frac{f''(0)}{2}.$$

As $\hat{f}_n^{\mathrm{S}}(0) = f(0) + \mathcal{O}_p(n^{-1/3})$, this implies that $\hat{B}_{21} n^{-1/3} = \mathcal{O}_p(n^{-2/9})$. Application of Theorem 3.1(i) yields that $\hat{f}_n^{\mathrm{A}}(0) = f(0) + \mathcal{O}_p(n^{-7/18})$. Sun and Woodroofe [12] also propose to use $\hat{f}_n^{\mathrm{P}}(\hat{\alpha}_n, 0)$ as an estimate of $f(0)$ in the case $k \geq 2$. They prove that in that case $n^{1/3}\{\hat{f}_n^{\mathrm{P}}(\hat{\alpha}_n, 0) - f(0)\} \to 0$ [see their Theorem 1(ii) on page 146].

We simulated 10,000 samples of sizes $n = 50$, 100, 200 and 10,000 from a standard half-normal distribution. For each sample the values were computed of $n^{1/3}\{\hat{f}_n^{\mathrm{S}}(0) - f(0)\}$, $n^{1/3}\{\hat{f}_n^{\mathrm{A}}(0) - f(0)\}$ and $n^{1/3}\{\hat{f}_n^{\mathrm{P}}(\hat{\alpha}_n, 0) - f(0)\}$. In Table 3 we list simulated values for the mean, variance and mean squared error of the three estimators. The simple estimator has the smallest variance, but as the sample size increases it becomes more biased. Nevertheless, it is stable for small sample sizes. The adaptive estimator becomes more biased with growing sample size, but with smaller MSE. The penalized MLE is most biased, also having a much larger variance than its simple and adaptive competitors.

TABLE 2
*Simulated mean, variance and mean squared error for both estimators at the half-normal distribution*

| | | $n$ | | | |
|---|---|---|---|---|---|
| | | **50** | **100** | **200** | **10,000** |
| $n^{2/5}\{\hat{f}_n^{\mathrm{S},2}(0) - f(0)\}$ | Mean | $-0.429$ | $-0.437$ | $-0.440$ | $-0.419$ |
| | Var | $0.371$ | $0.402$ | $0.440$ | $0.559$ |
| | MSE | $0.555$ | $0.592$ | $0.634$ | $0.735$ |
| $n^{2/5}\{\hat{f}_n^{\mathrm{A},2}(0) - f(0)\}$ | Mean | $-0.252$ | $-0.278$ | $-0.373$ | $-0.326$ |
| | Var | $0.459$ | $0.502$ | $0.549$ | $0.747$ |
| | MSE | $0.523$ | $0.579$ | $0.688$ | $0.853$ |

Finally, in Table 4 we list the true limiting values for the mean, variance and MSE, for all estimators at the exponential and half-normal distributions. The finite sample behavior of the simple estimators $\hat{f}_n^{\mathrm{S}}(0)$ (see Tables 1 and 3) and $\hat{f}_n^{\mathrm{S},2}(0)$ (see Table 2) reasonably matches the theoretical behavior. The adaptive estimators exhibit larger deviations from their theoretical values. This is probably explained by the fact that even for larger sample sizes, the estimation of the derivatives of $f$ in $B_{2k}$ still has a large influence.

One might prefer a scale-equivariant version of the above estimators. One possibility is $\hat{f}_n(X_{m:n})$, where $X_{m:n}$ denotes the $m$th order statistic. The sequence $m = m(n)$ should be chosen in such a way that $m(n) \to \infty$ and $m(n)/n \to 0$, for example, $m(n) = \lfloor an^{2/3} \rfloor$. In that case, one can show that $\hat{f}_n(X_{m:n})$ is asymptotically equivalent to $\hat{f}_n(af(0)^{-1}n^{-1/3})$. Its limiting distribution can be obtained from Theorem 3.1 and the AMSE optimal choice $a^*$ will depend on $f(0)$ and $f'(0)$. For this choice, $\hat{f}_n(a^* f(0)^{-1} n^{-1/3})$ has the same behavior as $\hat{f}_n(c_1^* B_{21} n^{-1/3})$. Another possibility is to estimate $f(0)$ by means of a numerical derivative of $F_n$,

$$\hat{f}_n^{\mathrm{D}}(0) = \frac{F_n(X_{m:n})}{X_{m:n}} = \frac{m/n}{X_{m:n}},$$

where $m = m(n)$ as above. It can be shown that $n^{1/3}\{\hat{f}_n^{\mathrm{D}}(0) - f(0)\}$ is asymptotically normal with mean $-|f'(0)|a/(2f(0))$ and variance $f(0)^2/a$. This implies that the minimal AMSE is a multiple of $(f(0)|f'(0)|)^{2/3}$, which also holds for $\hat{f}_n^{\mathrm{S}}(0)$ and $\hat{f}_n^{\mathrm{A}}(0)$ [see Theorem 3.1(ii) for the case $k = 1$]. Computer simulations show that the AMSE of $\hat{f}_n^{\mathrm{A}}(0)$ is always the smallest of the three.

TABLE 3
*Simulated mean, variance and mean squared error for the three estimators at the half-normal distribution*

| | | | | $n$ | |
|---|---|---|---|---|---|
| | | **50** | **100** | **200** | **10,000** |
| $n^{1/3}\{\hat{f}_n^{\mathrm{S}}(0) - f(0)\}$ | Mean | 0.012 | 0.058 | 0.104 | 0.269 |
| | Var | 0.320 | 0.317 | 0.316 | 0.296 |
| | MSE | 0.320 | 0.320 | 0.327 | 0.368 |
| $n^{1/3}\{\hat{f}_n^{\mathrm{A}}(0) - f(0)\}$ | Mean | 0.046 | 0.073 | 0.091 | 0.204 |
| | Var | 0.475 | 0.406 | 0.383 | 0.319 |
| | MSE | 0.477 | 0.412 | 0.391 | 0.361 |
| $n^{1/3}\{\hat{f}_n^{\mathrm{P}}(\hat{\alpha}_n, 0) - f(0)\}$ | Mean | 0.331 | 0.336 | 0.338 | 0.279 |
| | Var | 0.659 | 0.742 | 0.812 | 0.714 |
| | MSE | 0.768 | 0.855 | 0.926 | 0.792 |

TABLE 4
*Theoretical limiting mean, variance and mean squared error for the three estimators*

| | **Exponential** | | | **Half-normal** | | |
|---|---|---|---|---|---|---|
| **Estimator** | **Mean** | **Variance** | **MSE** | **Mean** | **Variance** | **MSE** |
| $n^{1/3}\{\hat{f}_n^{\mathrm{S}}(0) - f(0)\}$ | $-0.885$ | 0.805 | 1.591 | 0.336 | 0.316 | 0.429 |
| $n^{1/3}\{\hat{f}_n(c_1^* B_{21} n^{-1/3}) - f(0)\}$ | $-0.298$ | 1.043 | 1.131 | 0 | 0 | 0 |
| $n^{1/3}\{\hat{f}_n^{\mathrm{P}}(\hat{\alpha}_n, 0) - f(0)\}$ | $-0.349$ | 1.096 | 1.218 | 0 | 0 | 0 |
| $n^{2/5}\{\hat{f}_n^{\mathrm{S},2}(0) - f(0)\}$ | $-\infty$ | $\infty$ | $\infty$ | $-0.415$ | 0.670 | 0.842 |
| $n^{2/5}\{\hat{f}_n(c_2^* B_{22} n^{-1/5}) - f(0)\}$ | $-\infty$ | $\infty$ | $\infty$ | $-0.140$ | 0.718 | 0.737 |

## 6. Proofs.

PROOF OF LEMMA 2.1. Let $x_0 = \arg\max_{x \in K} f(x)$. If $x_0 = \infty$, there is nothing left to prove; therefore assume that $x_0 < \infty$.

(i) By definition of $x_0$ and the fact that $g$ is nonincreasing, for $x \geq x_0$, we must have $f(x) + g(x) \leq f(x_0) + g(x_0)$. Hence, we must have

$$\arg\max_{x \in K}\{f(x) + g(x)\} \leq x_0 = \arg\max_{x \in K} f(x).$$

This proves (i).

(ii) If $(C + x_0, \infty) \cap K = \varnothing$, the statement is trivially true, so only consider the case $(C + x_0, \infty) \cap K \neq \varnothing$. Then by definition $f(x) \leq f(x_0)$, for all $x \in (C + x_0, \infty) \cap K$, and by the property of $g$ we also have $g(x) \leq g(x_0)$, for $x \in (C + x_0, \infty) \cap K$. This implies $f(x) + g(x) \leq f(x_0) + g(x_0)$, for all $x \in (C + x_0, \infty) \cap K$. Hence, we must have

$$\arg\max_{x \in K}\{f(x) + g(x)\} \leq C + x_0 = C + \arg\max_{x \in K} f(x).$$

This proves the lemma. $\square$

PROOF OF LEMMA 3.1. Decompose the process $Z_{n1}$ as

$$
\begin{aligned}
(6.1) \quad Z_{n1}(x,t) &= n^{\alpha/2} W_n(F(tn^{-\alpha})) + n^{(1+\alpha)/2}\{F(tn^{-\alpha}) - f(0)tn^{-\alpha}\} \\
&\quad - xt - n^{\alpha/2} F(tn^{-\alpha}) W_n(1) + n^{\alpha/2} H_n(tn^{-\alpha}),
\end{aligned}
$$

where $H_n(t) = E_n(t) - B_n(F(t))$. By Brownian scaling, $n^{\alpha/2}W_n(F(tn^{-\alpha}))$ has the same distribution as the process $W(n^\alpha F(tn^{-\alpha}))$, and by uniform continuity of Brownian motion on compacta,

$$W(n^\alpha F(tn^{-\alpha})) - W(f(0)t) \to 0,$$

uniformly for $t$ in compact sets. Since $\alpha > 1/(2k+1)$ we have that

$$n^{(1+\alpha)/2}\{F(tn^{-\alpha}) - f(0)tn^{-\alpha}\} = n^{(1+\alpha)/2}\frac{f^{(k)}(\theta_t)}{(k+1)!}(tn^{-\alpha})^{k+1} \to 0,$$

uniformly for $t$ in compact sets. Because $n^{\alpha/2}F(tn^{-\alpha})W_n(1) = \mathcal{O}_p(n^{-\alpha/2})$, together with (2.2) this proves (i). In case (ii), where $\alpha = 1/(2k+1)$, the only difference is the behavior of the deterministic term

$$n^{(k+1)/(2k+1)}\{F(tn^{-1/(2k+1)}) - f(0)tn^{-1/(2k+1)}\} \to \frac{f^{(k)}(0)}{(k+1)!}t^{k+1},$$

uniformly for $t$ in compact sets. Similar to the proof of (i), using Brownian scaling and uniform continuity of Brownian motion on compacta this proves (ii).

For case (iii) the process $Z_{n2}$ can be written as

$$n^{b/2}\{W_n(F(cn^{-\alpha} + tn^{-b})) - W_n(F(cn^{-\alpha}))\}$$
$$+ n^{(b+1)/2}\{F(cn^{-\alpha} + tn^{-b}) - F(cn^{-\alpha}) - f(cn^{-\alpha})tn^{-b}\} - xt$$
$$- n^{b/2}\{F(cn^{-\alpha} + tn^{-b}) - F(cn^{-\alpha})\}W_n(1)$$
$$+ n^{b/2}H_n(cn^{-\alpha} + tn^{-b}) - n^{b/2}H_n(cn^{-\alpha}).$$

The process $n^{b/2}\{W_n(F(cn^{-\alpha} + tn^{-b})) - W_n(F(cn^{-\alpha}))\}$ has the same distribution as the process $W(n^b(F(cn^{-\alpha} + tn^{-b}) - F(cn^{-\alpha})))$, and by uniform continuity of Brownian motion on compacta,

$$W(n^b(F(cn^{-\alpha} + tn^{-b}) - F(cn^{-\alpha}))) - W(f(0)t) \to 0,$$

uniformly for $t$ in compact sets. Finally, for some $\theta_1 \in [cn^{-\alpha}, cn^{-\alpha} + tn^{-b}]$ and for some $\theta_2 \in [0, cn^{-\alpha} + tn^{-b}]$, it holds that

$$n^{(b+1)/2}\{F(cn^{-\alpha} + tn^{-b}) - F(cn^{-\alpha}) - f(cn^{-\alpha})tn^{-b}\}$$
$$= n^{(1-3b)/2}\frac{f'(\theta_1)}{2}t^2 = n^{(1-3b)/2}\frac{f^{(k)}(\theta_2)}{2(k-1)!}\theta_1^{k-1}t^2 \to \frac{f^{(k)}(0)}{2(k-1)!}c^{k-1}t^2,$$

uniformly for $t$ in compact sets. Since

$$n^{b/2}\{F(cn^{-\alpha} + tn^{-b}) - F(cn^{-\alpha})\}W_n(1) = \mathcal{O}_p(n^{-b/2}),$$

together with (2.2) this proves (iii). $\quad\square$

To verify condition (iii) of Theorem 2.1 we need that $F(c+t) - F(c) - f(c)t$ is suitably bounded. The next lemma guarantees that this is the case.

LEMMA 6.1. *Suppose that $f$ satisfies* (C2). *Then there exists a value $t_0 > 0$, such that $\inf |f^{(k)}| = \inf_{0 \le s \le t_0} |f^{(k)}(s)| > 0$. For any $0 \le c \le t_0/2$ we can bound $F(c+t) - F(c) - f(c)t$ by*

(i) $-\frac{\inf |f^{(k)}|}{(k+1)!} t^{k+1}$, *for $0 \le t \le t_0/2$,*

(ii) $-\frac{\inf |f^{(k)}|}{(k+1)!} (t_0/2)^k t$, *for $t > t_0/2$,*

(iii) $-\frac{\inf |f^{(k)}|}{2(k-1)!} (c/2)^{k-1} t^2$, *for $-c/2 < t < t_0/2$.*

*Furthermore, for small enough $c$ and for $-c < t < -c/2$,*

(iv) $F(c+t) - F(c) - f(c)t \le -C_1 c^{k+1}$, *where $C_1 > 0$ does not depend on $c$ and $t$.*

PROOF. The existence of $t_0 > 0$ follows directly from condition (C2). First note that if $f^{(k)}(0) \ne 0$, then we must have $f^{(k)}(0) < 0$, since otherwise $f^{(k-1)}$ is increasing in a neighborhood of zero, which implies that $f^{(k-2)}$ is increasing in a neighborhood of zero, and so on, which eventually would imply that $f$ is increasing in a neighborhood of zero. Therefore, under condition (C2) we must have $f^{(k)}(0) < 0$, which in turn implies that $f^{(i)}(s) < 0$ for $0 \le s \le t_0$ and $i = 1, 2, \ldots, k$. Hence, for $0 \le t \le t_0/2$, the inequality for $F(c+t) - F(c) - f(c)t$ is a direct consequence of a Taylor expansion, where all negative terms except for the last one are omitted.

For $t > t_0/2$, write

$$F(c+t) - F(c) - f(c)t$$
$$= F(c+t_0/2) - F(c) - f(c)t_0/2$$
$$\quad + (f(c+t_0/2) - f(c))(t - t_0/2)$$
$$\quad + F(c+t) - F(c+t_0/2) - f(c+t_0/2)(t - t_0/2),$$

where $F(c+t) - F(c+t_0/2) - f(c+t_0/2)(t - t_0/2) \le 0$, because $f$ is nonincreasing. By the same argument as above, $F(c+t_0/2) - F(c) - f(c)t_0/2 \le f^{(k)}(\theta_1)(t_0/2)^{k+1}/(k+1)!$ and $f(c+t_0/2) - f(c) \le f^{(k)}(\theta_2)(t_0/2)^k/k!$, for some $c < \theta_1, \theta_2 < c + t_0/2$. This implies that for $t > t_0/2$, we can bound $F(c+t) - F(c) - f(c)t$ from above by

$$-\frac{(t_0/2)^{k+1}}{(k+1)!} \inf |f^{(k)}| - \frac{(t_0/2)^k}{k!} \inf |f^{(k)}|(t - t_0/2)$$

$$\le -\frac{(t_0/2)^k}{(k+1)!} \inf |f^{(k)}|(t_0/2 + t - t_0/2)$$

$$= -\frac{(t_0/2)^k}{(k+1)!} \inf |f^{(k)}| t.$$

For $-c/2 < t < t_0/2$, first write $F(c + t) - F(c) - f(c)t = f'(\theta_4)t^2/2$, for $c/2 < \theta_4 < c + t_0/2$. By condition (C2), $f'(\theta_4) = f^{(k)}(\theta_5)\theta_4^{k-1}/(k-1)!$, for some $0 < \theta_5 < \theta_4$. This means that

$$F(c + t) - F(c) - f(c)t = \frac{\theta_4^{k-1}}{2(k-1)!}f^{(k)}(\theta_5)t^2 \leq -\frac{(c/2)^{k-1}}{2(k-1)!}\inf|f^{(k)}|t^2.$$

Finally, for $-c < t < -c/2$, first note that $f(c + t) - f(c) \geq 0$, so that $F(c + t) - F(c) - f(c)t$ is nondecreasing in $t$. Write

$$F(c + t) - F(c) - f(c)t$$
$$= \frac{f^{(k)}(\theta_6)}{(k+1)!}(c + t)^{k+1} - \frac{f^{(k)}(\theta_7)}{(k+1)!}c^{k+1} - \frac{f^{(k)}(\theta_8)}{k!}c^k t,$$

for $0 < \theta_6 < c + t$ and $0 < \theta_7, \theta_8 < c$. Because this expression is nondecreasing for $-c < t < -c/2$, and since $f^{(k)}(\theta_i) - f^{(k)}(0) = o(1)$, for $i = 6, 7, 8$, uniformly in $-c < t < -c/2$, we conclude that

$$F(c + t) - F(c) - f(c)t \leq \frac{f^{(k)}(0)}{(k+1)!}c^{k+1}\left(\frac{1}{2^{k+1}} - 1 + \frac{k+1}{2}\right)(1 + o(1))$$

as $c \downarrow 0$. Since $f^{(k)}(0) < 0$, this proves the lemma. $\square$

PROOF OF LEMMA 3.2. (i) Decompose $Z_{n1}$ as in (6.1). Let $0 < \varepsilon < x$ and define

$$X_{n1}(t) = n^{\alpha/2}H_n(tn^{-\alpha}) - \varepsilon t/2,$$

where $H_n(t) = E_n(t) - B_n(F(t))$. Next, consider the event

$$(6.2) \quad A_{n1} = \{X_{n1}(s) \geq X_{n1}(t), \text{ for all } s, t \geq 0, \text{ such that } t - s \geq \delta_n\}.$$

Then with $\delta_n = n^{-(1-\alpha)/2}(\log n)^2$, by using (2.2) we have that

$$P(A_{n1}) \geq P\left\{\sup_{t \in [0, \infty)}|H_n(t)| \leq \frac{\varepsilon}{4}n^{-1/2}(\log n)^2\right\} \to 1.$$

Also define the process $X_{n2}(t) = -n^{\alpha/2}F(tn^{-\alpha})W_n(1) - \varepsilon t/2$, and consider the event

$$(6.3) \qquad A_{n2} = \{X_{n2}(s) \geq X_{n2}(t), \text{ for all } 0 \leq s \leq t < \infty\}.$$

Then, since every sample path of the process $X_{n2}$ is differentiable, we have

$$P(A_{n2}) \geq P\left\{-f(tn^{-\alpha})W_n(1) - \frac{\varepsilon}{2}n^{\alpha/2} \leq 0, \text{ for all } t \in [0, \infty)\right\} \to 1.$$

Hence, if $A_n = A_{n1} \cap A_{n2}$, then $P(A_n) \to 1$. Since for any $\eta > 0$,

$$P\left\{\underset{t \in [0, \infty)}{\arg\max} Z_{n1}(t)\mathbb{1}_{A_n^c} > \eta\right\} \leq P(A_n^c) \to 0,$$

we conclude that $(\arg\max_t Z_{n1}(t))\mathbb{1}_{A_n^c} = \mathcal{O}_p(1)$. This means that we only have to consider $(\arg\max_t Z_{n1}(t))\mathbb{1}_{A_n}$. From Lemma 2.1 we have

$$(6.4) \qquad \left(\arg\max_{t\in[0,\infty)} Z_{n1}(t)\right)\mathbb{1}_{A_n} \le \arg\max_{t\in[0,\infty)} S_{n1}(t) + \delta_n,$$

where

$$S_{n1}(t) = n^{\alpha/2} W_n(F(tn^{-\alpha})) - (x-\varepsilon)t + n^{(1+\alpha)/2}(F(tn^{-\alpha}) - f(0)tn^{-\alpha}).$$

Since $F(tn^{-\alpha}) - f(0)tn^{-\alpha}$ is nonincreasing for $t \ge 0$, according to Lemma 2.1,

$$
\begin{aligned}
(6.5) \qquad \arg\max_{t\in[0,\infty)} S_{n1}(t) &\le \arg\max_{t\in[0,\infty)}\{n^{\alpha/2} W_n(F(tn^{-\alpha})) - (x-\varepsilon)t\} \\
&\le \sup\{t \ge 0 : n^{\alpha/2} W_n(F(tn^{-\alpha})) - (x-\varepsilon)t \ge 0\}.
\end{aligned}
$$

By change of variables $u = G(t) = n^\alpha F(tn^{-\alpha})$, and using that for $u \in [0, n^\alpha]$,

$$(6.6) \qquad \frac{u}{f(0)} \le G^{-1}(u) \le \frac{u}{f(F^{-1}(un^{-\alpha}))},$$

we find that the right-hand side of (6.5) is bounded by

$$G^{-1}\left(\sup\left\{u \ge 0 : n^{\alpha/2} W_n(un^{-\alpha}) - \frac{x-\varepsilon}{f(0)}u \ge 0\right\}\right).$$

By Brownian scaling (2.3),

$$\sup\left\{u \ge 0 : n^{\alpha/2} W_n(un^{-\alpha}) - \frac{x-\varepsilon}{f(0)}u \ge 0\right\} \stackrel{d}{=} \frac{f(0)^2}{(x-\varepsilon)^2}\sup\{u \ge 0 : W(u) - u \ge 0\},$$

which is of order $\mathcal{O}_p(1)$. The latter can be seen, for instance, from (3.3). Because $\delta_n = n^{-(1-\alpha)/2}(\log n)^2 = o(1)$, together with (6.4), (6.5) and (6.6), it follows that

$$
\begin{aligned}
0 \le \arg\max_{t\in[0,\infty)} Z_{n1}(t) &\le \left(\arg\max_{t\in[0,\infty)} Z_{n1}(t)\right)\mathbb{1}_{A_n} + \mathcal{O}_p(1) \\
&\le \frac{\mathcal{O}_p(1)}{f(F^{-1}(\mathcal{O}_p(n^{-\alpha})))} + \mathcal{O}_p(1),
\end{aligned}
$$

which proves (i).

(ii) In this case $\alpha = 1/(2k+1)$, so that the argument up to (6.4) is the same. Let $\varepsilon > 0$ and $A_n = A_{n1} \cap A_{n2}$, where $A_{n1}$ is as defined in (6.2) with $\delta_n = n^{-k/(k+1)}(\log n)^2$ and $A_{n2}$ is as defined in (6.3). We now find that

$$
\begin{aligned}
(6.7) \qquad \left(\arg\max_{t\in[0,\infty)} Z_{n1}(t)\right)\mathbb{1}_{A_n} &\le \arg\max_{t\in[0,\infty)} S_{n1}(t) + \delta_n \\
&\le \sup\{t \ge 0 : S_{n1}(t) \ge 0\} + \delta_n.
\end{aligned}
$$

Let $t_0$ be the value from Lemma 6.1 and consider the event

$$D_{n1} = \{n^{-\alpha} \sup\{t \geq 0 : S_{n1}(t) \geq 0\} \leq t_0/2\}.$$

If $S_{n1}(t) \geq 0$, then according to Lemma 6.1(ii), for $tn^{-\alpha} > t_0/2$ and $n$ sufficiently large, we find that

$$0 \leq n^{\alpha/2} W_n(F(tn^{-\alpha})) - (x - \varepsilon)t + n^{(1+\alpha)/2}(F(tn^{-\alpha}) - f(0)tn^{-\alpha})$$

$$\leq n^{\alpha/2} \sup_{0 \leq u \leq 1} |W_n(u)| - (x - \varepsilon)t - n^{(1-\alpha)/2}\frac{(t_0/2)^k}{(k+1)!}\inf|f^{(k)}|t$$

$$\leq n^{\alpha/2} \sup_{0 \leq u \leq 1} |W_n(u)| - n^{(1-\alpha)/2}C_1 t\left(1 + \frac{x - \varepsilon}{n^{(1-\alpha)/2}C_1}\right)$$

$$\leq n^{\alpha/2}\left\{\sup_{0 \leq u \leq 1} |W_n(u)| - C_1 n^{1/2}t_0/4\right\},$$

where $C_1 = \inf|f^{(k)}|(t_0/2)^k/(k+1)!$. Therefore

$$P(D_{n1}^c) \leq P\left(\sup_{0 \leq u \leq 1}|W(u)| \geq C_1 n^{1/2}t_0/4\right) \to 0.$$

This means we can restrict ourselves to the event $A_n \cap D_{n1}$, so that by reasoning analogous to that before, from (6.7) we get

$$\left(\operatorname*{arg\,max}_{t \in [0, \infty)} Z_{n1}(t)\right)\mathbb{1}_{A_n \cap D_{n1}} \leq \sup\{t \geq 0 : S_{n1}(t) \geq 0\}\mathbb{1}_{D_{n1}} + \delta_n$$

$$\leq \sup\{0 \leq t \leq n^\alpha t_0/2 : S_{n1}(t) \geq 0\} + \delta_n.$$

According to Lemma 6.1(i), for $0 \leq tn^{-\alpha} \leq t_0/2$ and using that $\alpha = 1/(2k + 1)$, we get

$$n^{(1+\alpha)/2}(F(tn^{-\alpha}) - f(0)tn^{-\alpha}) \leq -\frac{\inf|f^{(k)}|}{(k+1)!}t^{k+1},$$

so that

$$0 \leq \left(\operatorname*{arg\,max}_{t \in [0, \infty)} Z_{n1}(t)\right)\mathbb{1}_{A_n \cap D_{n1}}$$

$$(6.8) \qquad \leq \sup\left\{0 \leq t \leq n^\alpha t_0/2 : n^{\alpha/2} W_n(F(tn^{-\alpha}))\right.$$

$$\left. - (x - \varepsilon)t - \frac{\inf|f^{(k)}|}{(k+1)!}t^{k+1} \geq 0\right\} + \delta_n.$$

Next, distinguish between

(A)  $-(x - \varepsilon)t - \inf|f^{(k)}|t^{k+1}/(2(k+1)!) \geq 0,$

(B) $-(x-\varepsilon)t - \inf|f^{(k)}|t^{k+1}/(2(k+1)!) < 0$.

Since $t \geq 0$, case (A) can only occur when $x - \varepsilon < 0$, in which case we have $0 \leq t \leq (2(k+1)!(\varepsilon-x)/\inf|f^{(k)}|)^{1/k}$, which is of order $\mathcal{O}(1)$. In case (B), it follows that

$$n^{\alpha/2}W_n(F(tn^{-\alpha})) - \frac{\inf|f^{(k)}|}{2(k+1)!}t^{k+1} \geq 0.$$

We conclude from (6.8) that

$$0 \leq \left(\underset{t\in[0,\infty)}{\arg\max}\, Z_{n1}(t)\right)\mathbb{1}_{A_n\cap D_{n1}}$$

$$\leq \sup\left\{0 \leq t \leq n^{\alpha}t_0/2 : n^{\alpha/2}W_n(F(tn^{-\alpha})) - \frac{\inf|f^{(k)}|}{2(k+1)!}t^{k+1} \geq 0\right\}$$

(6.9) $$\quad + \mathcal{O}_p(1) + \delta_n$$

$$\leq \sup\left\{t \in [0,\infty) : n^{\alpha/2}W_n(F(tn^{-\alpha})) - \frac{\inf|f^{(k)}|}{2(k+1)!}t^{k+1} \geq 0\right\}$$

$$\quad + \mathcal{O}_p(1).$$

Similar to the proof of (i), by change of variables $u = G(t) = n^{\alpha}F(tn^{-\alpha})$ and using (6.6) with $\alpha = 1/(2k+1)$, we find that the arg max on the right-hand side of (6.9) is bounded from above by

$$G^{-1}\left(\sup\left\{u \in [0,\infty) : n^{\alpha/2}W_n(un^{-\alpha}) - \frac{\inf|f^{(k)}|u^{k+1}}{2(k+1)!f(0)^{k+1}} \geq 0\right\}\right) + \mathcal{O}_p(1).$$

By Brownian scaling (2.3), we obtain that the supremum in the first term has the same distribution as

$$\left(\frac{2(k+1)!f(0)^{k+1}}{\inf|f^{(k)}|}\right)^{2/(2k+1)}\sup\{u \geq 0 : W(u) - u^{k+1} \geq 0\}.$$

Again by using (3.3), this is of order $\mathcal{O}_p(1)$. Similar to the proof of (i), from (6.6) and (6.9) we find that

$$0 \leq \underset{t\in[0,\infty)}{\arg\max}\, Z_{n1}(t) \leq \left(\underset{t\in[0,\infty)}{\arg\max}\, Z_{n1}(t)\right)\mathbb{1}_{A_n\cap D_{n1}} + \mathcal{O}_p(1)$$

$$\leq \frac{\mathcal{O}_p(1)}{f(F^{-1}(\mathcal{O}_p(n^{-\alpha})))} + \mathcal{O}_p(1),$$

which proves (ii).

(iii) Decompose $Z_{n2}$ as in the proof of Lemma 3.1. Let $\varepsilon > 0$ and $A_n = A_{n1} \cap A_{n2}$, with $A_{n1}$ defined similarly to (6.2) with $\delta_n = n^{-(1-b)/2}(\log n)^2$, where $b$ is the same as in Lemma 3.1, and $A_{n2}$ is defined similarly to (6.3).

By the same argument as in the proof of (i) and (ii), it suffices to consider $(\arg\max_t Z_{n2}(t))\mathbb{1}_{A_n}$. We find

$$\left(\underset{t\in[-cn^{b-\alpha},\infty)}{\arg\max} Z_{n2}(t)\right)\mathbb{1}_{A_n} \le \underset{t\in[-cn^{b-\alpha},\infty)}{\arg\max} M_{n2}(t) + \delta_n$$

$$\le \sup\{t \ge 0 : M_{n2}(t) \ge 0\} + \delta_n,$$

where $M_{n2}(t)$ has the same distribution as

$$S_{n2}(t) = n^{b/2}W(F(cn^{-\alpha} + tn^{-b}) - F(cn^{-\alpha}))$$

$$+ n^{(b+1)/2}(F(cn^{-\alpha} + tn^{-b}) - F(cn^{-\alpha}) - f(cn^{-\alpha})tn^{-b})$$

$$- (x - \varepsilon)t.$$

As in the proof of (ii), consider $D_{n2} = \{n^{-b}\sup\{t \ge 0 : S_{n2}(t) \ge 0\} \le t_0/2\}$, where $t_0$ is the value from Lemma 6.1. By the same reasoning as used in the proof of (ii), it again follows from Lemma 6.1(ii) that $P(D_{n2}^c) \to 0$, so we only have to consider $\sup\{t \ge 0 : S_{n2}(t) \ge 0\}\mathbb{1}_{D_{n2}}$. Hence, similar to the proof of (ii) we get

$$\sup\{t \ge 0 : S_{n2}(t) \ge 0\}\mathbb{1}_{D_{n2}} \le \sup\{0 \le t \le n^b t_0/2 : S_{n2}(t) \ge 0\}.$$

Since $b > 1/(2k + 1)$, for $k \ge 2$, we cannot proceed as in the proof of (ii) by using Lemma 6.1(i) to bound the drift term. However, according to Lemma 6.1(iii), for $0 \le t \le n^b t_0/2$,

$$n^{(b+1)/2}(F(cn^{-\alpha} + tn^{-b}) - F(cn^{-\alpha}) - f(cn^{-\alpha})tn^{-b}) \le -\frac{\inf|f^{(k)}|}{2^k(k-1)!}t^2,$$

so that $\sup\{0 \le t \le n^b t_0/2 : S_{n2}(t) \ge 0\}$ is bounded from above by

$$\sup\left\{0 \le t \le n^b t_0/2 : n^{b/2}W(F(cn^{-\alpha} + tn^{-b}) - F(cn^{-\alpha}))\right.$$

$$\left. - (x - \varepsilon)t - \frac{\inf|f^{(k)}|}{2^k(k-1)!}t^2 \ge 0\right\}.$$

Similarly to (6.9), we conclude that $\sup\{t \ge 0 : S_{n2}(t) \ge 0\}\mathbb{1}_{D_{n2}}$ is bounded from above by

$$(6.10) \qquad \sup\left\{t \ge 0 : n^{b/2}W_n(F(cn^{-\alpha} + tn^{-b}) - F(cn^{-\alpha}))\right.$$

$$\left. - \frac{\inf|f^{(k)}|}{2^{k+1}(k-1)!}t^2 \ge 0\right\} + \mathcal{O}_p(1).$$

Next, change variables $u = G(t) = n^b(F(cn^{-\alpha} + tn^{-b}) - F(cn^{-\alpha}))$. Then for any $u \in [0, n^b(1 - F(cn^{-\alpha}))]$, it follows that

$$(6.11) \qquad \frac{u}{f(0)} \le G^{-1}(u) \le \frac{u}{f(F^{-1}(un^{-b} + F(cn^{-\alpha})))},$$

so that (6.10) is bounded from above by

$$G^{-1}\left(\sup\left\{u \geq 0 : n^{b/2}W(un^{-b}) - \frac{\inf|f^{(k)}|}{2^{k+1}(k-1)!f(0)^2}u^2 \geq 0\right\}\right) + \mathcal{O}_p(1).$$

As in the proof of (ii), by Brownian scaling (2.3) together with (6.11), we find that

$$\begin{aligned}
(6.12) \quad \underset{t \in [-cn^{b-\alpha}, \infty)}{\arg\max} Z_{n2}(t) &\leq \left(\underset{t \in [-cn^{b-\alpha}, \infty)}{\arg\max} Z_{n2}(t)\right)\mathbb{1}_{A_n \cap D_{n2}} + \mathcal{O}_p(1) \\
&\leq \frac{\mathcal{O}_p(1)}{f(F^{-1}(\mathcal{O}_p(n^{-b}) + F(cn^{-\alpha})))} + \mathcal{O}_p(1) \\
&= \mathcal{O}_p(1).
\end{aligned}$$

To obtain a lower bound for the left-hand side of (6.12), first note that

$$(6.13) \quad \underset{t \in [-cn^{b-\alpha}, \infty)}{\arg\max} Z_{n2}(t) \geq \underset{t \in [-cn^{b-\alpha}, 0]}{\arg\max} Z_{n2}(t) = -\underset{t \in [0, cn^{b-\alpha}]}{\arg\max} Z_{n2}(-t).$$

From here the argument runs along the same lines as for the upper bound. Let $\varepsilon > 0$ and, similarly to (6.2) and (6.3), define the events $A_{n1}$ and $A_{n2}$ with

$$X_{n1}(t) = n^{b/2}H_n(cn^{-\alpha} - tn^{-b}) - \varepsilon t/2,$$

$$X_{n2}(t) = -n^{b/2}F(cn^{-\alpha} - tn^{-b})W_n(1) - \varepsilon t/2.$$

With $A_n = A_{n1} \cap A_{n2}$, as before we get $(\arg\max_t Z_{n2}(-t))\mathbb{1}_{A_n^c} = \mathcal{O}_p(1)$ and

$$\left(\arg\max_t Z_{n2}(-t)\right)\mathbb{1}_{A_n} \leq \underset{t \in [0, cn^{b-\alpha})}{\arg\max} M_{n3}(t) + \delta_n,$$

where $M_{n3}(t)$ has the same distribution as

$$\begin{aligned}
S_{n3}(t) &= n^{b/2}W(F(cn^{-\alpha} - tn^{-b}) - F(cn^{-\alpha})) \\
&\quad + n^{(b+1)/2}(F(cn^{-\alpha} - tn^{-b}) - F(cn^{-\alpha}) + f(cn^{-\alpha})tn^{-b}) \\
&\quad + (x + \varepsilon)t \\
&\leq n^{b/2}\sup\{|W(u)| : 0 \leq u \leq f(0)cn^{-\alpha}\} \\
&\quad + n^{(b+1)/2}(F(cn^{-\alpha} - tn^{-b}) - F(cn^{-\alpha}) + f(cn^{-\alpha})tn^{-b}) \\
&\quad + (x + \varepsilon)t.
\end{aligned}$$

Consider $D_{n3} = \{n^{-b}\sup\{0 \leq t \leq cn^{b-\alpha} : S_{n3}(t) \geq 0\} \leq cn^{-\alpha}/2\}$, and note that by Brownian scaling $\sup\{|W(u)| : 0 \leq u \leq f(0)cn^{-\alpha}\}$ has the same distribution as $n^{-\alpha/2}\sup\{|W(u)| : 0 \leq u \leq cf(0)\}$. Reasoning as in the proof of

(ii), using Lemma 6.1(iv), we obtain that for $cn^{-\alpha}/2 \le n^{-b}t \le cn^{-\alpha}$ and $n$ sufficiently large,

$$0 \le n^{(b-\alpha)/2} \sup_{0 \le u \le cf(0)} |W(u)|$$

$$+ n^{(b+1)/2}(F(cn^{-\alpha} - tn^{-b}) - F(cn^{-\alpha}) + f(cn^{-\alpha})tn^{-b}) + (x+\varepsilon)t$$

$$\le n^{(b-\alpha)/2}\Bigg( \sup_{0 \le u \le cf(0)} |W(u)|$$

$$- C_1 n^{(1-(2k+1)\alpha)/2}\Big(1 + \frac{x+\varepsilon}{C_1 n^{(b+1)/2 - (k+1)\alpha}}\Big)\Bigg)$$

$$\le n^{(b-\alpha)/2}\Bigg( \sup_{0 \le u \le cf(0)} |W(u)| - \frac{C_1}{2} n^{(1-(2k+1)\alpha)/2}\Bigg).$$

Therefore, $P(D_{n3}^c) \to 0$, so we only have to consider $(\arg\max_t S_{n3}(t))\mathbb{1}_{D_{n3}}$. Hence, similar to the proof of (ii), we get

$$\underset{t \in [0, cn^{b-\alpha})}{\arg\max} \, S_{n3}(t)\mathbb{1}_{D_{n3}} + \delta_n \le \sup\{0 \le t \le cn^{b-\alpha}/2 : S_{n3}(t) \ge 0\} + \delta_n.$$

According to Lemma 6.1(iii), for $0 \le tn^{-b} \le cn^{-\alpha}/2$ we have

$$(6.14) \quad \begin{aligned} & n^{(b+1)/2}(F(cn^{-\alpha} - tn^{-b}) - F(cn^{-\alpha}) + f(cn^{-\alpha})tn^{-b}) \\ & \le -\frac{\inf |f^{(k)}|}{2^k(k-1)!}t^2. \end{aligned}$$

Similar to (ii), separate cases and obtain that $\arg\max_{t \in [0, cn^{b-\alpha})} S_{n3}(t)\mathbb{1}_{D_{n3}} + \delta_n$ is bounded from above by

$$\sup\Bigg\{ 0 \le t \le cn^{b-\alpha}/2 : n^{b/2}W(F(cn^{-\alpha} - tn^{-b}) - F(cn^{-\alpha}))$$

$$- \frac{\inf |f^{(k)}|}{2^{k+1}(k-1)!}t^2 \ge 0 \Bigg\} + \mathcal{O}_p(1).$$

After change of variables $u = G(t) = n^b(F(cn^{-\alpha} - tn^{-b}) - F(cn^{-\alpha}))$, and using that $u \in [-n^b F(cn^{-\alpha}), 0]$, one has

$$-\frac{u}{f(0)} \le G^{-1}(u) \le -\frac{u}{f(cn^{-\alpha})}.$$

We now find that

$$\underset{t \in [0, cn^{b-\alpha})}{\arg\max} \, S_{n3}(t) + \delta_n$$

$$\le \frac{1}{f(cn^{-\alpha})} \sup\Bigg\{ u \le 0 : W_n(u) - \frac{\inf |f^{(k)}|}{2^{k+1}(k-1)!f(0)^2}u^2 \ge 0 \Bigg\} + \mathcal{O}_p(1).$$

As above, by Brownian scaling (2.3) together with (6.13), it follows that

$$\underset{t \in [-cn^{b-\alpha}, \infty)}{\arg\max} \; Z_{n2}(t) \geq \frac{\mathcal{O}_p(1)}{f(cn^{-\alpha})} + \mathcal{O}_p(1) = \mathcal{O}_p(1).$$

Together with (6.12) this proves the lemma. $\square$

**Acknowledgments.** The authors wish to thank two unknown referees for their fruitful suggestions that have lead to a substantial improvement of the original manuscript. We thank the first referee for his/her careful reading of the manuscript and for providing a neat closed expression of the limit distribution in Theorem 3.1(ii) and suggesting to incorporate vanishing derivatives at zero. We thank the second referee for mentioning the scale-equivariant alternatives and suggesting to incorporate the constant $c$ in Theorem 3.1, which inspired us to take a closer look at the problem of estimating $f(0)$.

ING FINANCIAL MARKETS
FOPPINGADREEF 7
P.O. BOX 1800
1000 BV AMSTERDAM
THE NETHERLANDS
E-MAIL: vladimir.kulikov@ingbank.com

DELFT INSTITUTE OF APPLIED MATHEMATICS
FACULTY OF EEMCS
MEKELWEG 4
2628 CD DELFT
THE NETHERLANDS
E-MAIL: h.p.lopuhaa@ewi.tudelft.nl